\documentclass[a4paper,12pt]{article}
\usepackage[T1]{fontenc}
\usepackage[latin1]{inputenc}
\usepackage[english]{babel}
\usepackage{hyperref}
\usepackage[dvips]{graphics}
\usepackage{amsmath}
\usepackage{amssymb}
\usepackage{amsopn}
\usepackage{amsbsy}
\usepackage{amstext}
\usepackage{latexsym}
\usepackage{amsfonts}
\usepackage{array}
\usepackage{epsfig}
\usepackage{color}
\usepackage{mathrsfs} 
\title{A study of a curious arithmetic function}
\author{{\sc Bakir FARHI} \\ {\tt bakir.farhi@gmail.com}}
\date{}

\newtheorem{thm}{Theorem}[section]

\newtheorem{coll}[thm]{Corollary}
\newtheorem{conj}[thm]{Conjecture}

\let\epsilon=\varepsilon

\def\c{{\rm c}}
\def\lcm{{\rm lcm}}
\def\odd{{\rm Odd}}

\def\EMdash{\leavevmode\hbox to 7.5mm{\vrule height .63ex depth -.59ex
    width 5.4mm\hfill}}

\begin{document}
\maketitle \vspace{-6cm}
\begin{flushleft}
To appear
\end{flushleft}~\vspace{4cm}

\begin{abstract}
In this note, we study the arithmetic function $f : \mathbb{Z}_+^*
\rightarrow \mathbb{Q}_+^*$ defined by $f(2^k \ell) = \ell^{1 -
k}$ ($\forall k , \ell \in \mathbb{N}$, $\ell$ odd). We show
several important properties about that function and then we use
them to obtain some curious results involving the $2$-adic
valuation.
\end{abstract}
\vspace{1mm}

\section{Introduction and notations}~

The purpose of this paper is to study the arithmetic function $f :
\mathbb{Z}_+^* \rightarrow \mathbb{Q}_+^*$ defined by:
$$f(2^k \ell) = \ell^{1-k} ~~~~ (\forall k , \ell \in \mathbb{N} , \ell ~\text{odd}) .$$
We have for example $f(1) = 1 , f(2) = 1 , f(3) = 3 , f(12) =
\frac{1}{3} , f(40) = \frac{1}{25} , \dots$. So it is clear that
$f(n)$ is not always an integer. However, we will show in what
follows that $f$ satisfies among others the property that the
product of the $f(r)$'s ($1 \leq r \leq n$) is always an integer
and it is a multiple of all odd prime number not exceeding $n$.
Further, we exploit the properties of $f$ to establish some
curious properties concerning the $2$-adic valuation.

The study of $f$ requires to introduce the two auxiliary
arithmetic functions $g : \mathbb{Q}_+^* \rightarrow
\mathbb{Z}_+^*$ and $h : \mathbb{Z}_+^* \rightarrow
\mathbb{Q}_+^*$, defined by:
\begin{equation}\label{eq1}
g(x) := \begin{cases} x & \text{if $x \in \mathbb{N}$} \\
1 & \text{else}
\end{cases}~~~~~~~~ (\forall x \in \mathbb{Q}_+^*)
\end{equation}
\begin{equation}\label{eq2}
h(r) := \frac{r}{g(\frac{r}{2}) g(\frac{r}{4}) g(\frac{r}{8})
\cdots} ~~~~~~~~ (\forall r \in \mathbb{Z}_+^*)
\end{equation}
Remark that the product in the denominator of the right-hand side
of (\ref{eq2}) is actually finite because $g(\frac{r}{2^i}) = 1$
for any sufficiently large $i$; so $h$ is
well-defined.~\vspace{2mm}

\noindent{\bf Some notations and terminologies.} Throughout this
paper, we let $\mathbb{N}^*$ denote the set $\mathbb{N} \setminus
\{0\}$ of positive integers. For a given prime number $p$, we let
$v_p$ denote the usual $p$-adic valuation. We define the {\it odd
part} of a positive rational number $\alpha$ as the positive
rational number, denoted $\odd(\alpha)$, so that we have $\alpha =
2^{v_2(\alpha)} \cdot \odd(\alpha)$. Finally, we denote by
$\lfloor . \rfloor$ the integer-part function and we often use in
this paper the following elementary well-known property of that
function:
$$\forall a , b \in \mathbb{N}^* , \forall x \in \mathbb{R} :~~~~ \left\lfloor\frac{\left\lfloor\frac{x}{a}\right\rfloor}{b}\right\rfloor
= \left\lfloor\frac{x}{a b}\right\rfloor .$$

\section{Results and proofs}

\begin{thm}\label{t1}
Let $n$ be a positive integer. Then the product $\displaystyle
\prod_{r = 1}^{n} f(r)$ is an integer.
\end{thm}
{\bf Proof.} For a given $r \in \mathbb{N}^*$, let us write $f(r)$
in terms of $h(r)$. By writing $r$ in the form $r = 2^k \ell$ ($k
, \ell \in \mathbb{N}$, $\ell$ odd), we have by the definition of
$g$:
$$g\left(\frac{r}{2}\right) g\left(\frac{r}{4}\right) g\left(\frac{r}{8}\right) \cdots = \left(2^{k - 1} \ell\right) (2^{k - 2} \ell) \times \cdots \times (2^0 \ell) =
2^{\frac{k (k - 1)}{2}} \ell^k .$$ So, it follows that:
$$h(r) := \frac{r}{g(\frac{r}{2}) g(\frac{r}{4}) g(\frac{r}{8}) \cdots} = \frac{2^k \ell}{2^{\frac{k (k - 1)}{2}} \ell^k}
= 2^{\frac{k (3 - k)}{2}} \ell^{1 - k} = 2^{\frac{k (3 - k)}{2}}
f(r) .$$ Hence
\begin{equation}\label{eq3}
f(r) = 2^{\frac{v_2(r) (v_2(r) - 3)}{2}} h(r) .
\end{equation}
Using (\ref{eq3}), we get for all $n \in \mathbb{N}^*$:
\begin{equation}\label{eq4}
\prod_{r = 1}^{n} f(r) ~=~ 2^{\sum_{r = 1}^{n}\frac{v_2(r) (v_2(r)
- 3)}{2}} \prod_{r = 1}^{n} h(r) .
\end{equation}
By taking the odd part of each of the two hand-sides of this last
identity, we obtain:
\begin{equation}\label{eq5}
\prod_{r = 1}^{n} f(r) = \odd\left(\prod_{r = 1}^{n} h(r)\right)
~~~~ (\forall n \in \mathbb{N}^*) .
\end{equation}
So, to confirm the statement of the theorem, it suffices to prove
that the product $\prod_{r = 1}^{n} h(r)$ is an integer for any $n
\in \mathbb{N}^*$. To do so, we lean on the following sample
property of $g$:
$$g\left(\frac{1}{a}\right) g\left(\frac{2}{a}\right) \cdots g\left(\frac{r}{a}\right) = \left\lfloor \frac{r}{a} \right\rfloor ! ~~~~ (\forall r , a
\in \mathbb{N}^*) .$$ Using this, we have:
\begin{eqnarray*}
\prod_{r = 1}^{n} h(r) & = & \prod_{r = 1}^{n}
\frac{r}{g\left(\frac{r}{2}\right) g\left(\frac{r}{4}\right)
g\left(\frac{r}{8}\right) \cdots} \\
& = & \frac{n!}{\displaystyle\prod_{r = 1}^{n}
g\left(\frac{r}{2}\right) \cdot \displaystyle\prod_{r = 1}^{n}
g\left(\frac{r}{4}\right) \cdot \prod_{r = 1}^{n}
g\left(\frac{r}{8}\right) \cdots} \\
& = & \frac{n!}{\lfloor\frac{n}{2}\rfloor !
\lfloor\frac{n}{4}\rfloor ! \lfloor\frac{n}{8}\rfloor ! \cdots} .
\end{eqnarray*}
Hence
\begin{equation}\label{eq6}
\prod_{r = 1}^{n} h(r) = \frac{n!}{\lfloor\frac{n}{2}\rfloor !
\lfloor\frac{n}{4}\rfloor ! \lfloor\frac{n}{8}\rfloor ! \cdots}
\end{equation}
(Remark that the product in the denominator of the right-hand side
of (\ref{eq6}) is actually finite because
$\lfloor\frac{n}{2^i}\rfloor = 0$ for any sufficiently large
$i$).\\
Now, since $\lfloor\frac{n}{2}\rfloor + \lfloor\frac{n}{4}\rfloor
+ \lfloor\frac{n}{8}\rfloor + \dots \leq \frac{n}{2} + \frac{n}{4}
+ \frac{n}{8} + \dots = n$ then
$\frac{n!}{\lfloor\frac{n}{2}\rfloor ! \lfloor\frac{n}{4}\rfloor !
\lfloor\frac{n}{8}\rfloor ! \cdots}$ is a multiple of the
multinomial coefficient $\binom{\lfloor\frac{n}{2}\rfloor +
\lfloor\frac{n}{4}\rfloor + \lfloor\frac{n}{8}\rfloor +
\dots}{\lfloor\frac{n}{2}\rfloor ~ \lfloor\frac{n}{4}\rfloor ~
\lfloor\frac{n}{8}\rfloor ~ \dots}$ which is an integer.
Consequently $\frac{n!}{\lfloor\frac{n}{2}\rfloor !
\lfloor\frac{n}{4}\rfloor ! \lfloor\frac{n}{8}\rfloor ! \cdots}$
is an integer, which completes this proof.\hfill$\blacksquare$

\begin{thm}\label{t2}
Let $n$ be a positive integer. Then $\displaystyle \prod_{r =
1}^{n} f(r)$ is a multiple of $\odd(\lcm(1 , 2 , \dots , n))$.\\
In particular, $\displaystyle \prod_{r = 1}^{n} f(r)$ is a
multiple of all odd prime number not exceeding $n$.
\end{thm}

\noindent{\bf Proof.} According to the relations (\ref{eq5}) and
(\ref{eq6}) obtained during the proof of Theorem \ref{t1}, it
suffices to show that $\frac{n!}{\lfloor\frac{n}{2}\rfloor !
\lfloor\frac{n}{4}\rfloor ! \lfloor\frac{n}{8}\rfloor ! \cdots}$
is a multiple of $\lcm(1 , 2 , \dots , n)$. Equivalently, it
suffices to prove that for all prime number $p$, we have:
\begin{equation}\label{eq7}
v_p\left(\frac{n!}{\lfloor\frac{n}{2}\rfloor !
\lfloor\frac{n}{4}\rfloor ! \lfloor\frac{n}{8}\rfloor !
\cdots}\right) \geq \alpha_p ,
\end{equation}
where $\alpha_p$ is the $p$-adic valuation of $\lcm(1 , 2 , \dots
, n)$, that is the greatest power of $p$ not exceeding $n$. Let us
show (\ref{eq7}) for a given arbitrary prime number $p$. Using
Legendre's formula (see e.g., \cite{hw}), we have:
\begin{eqnarray}
v_p\left(\frac{n!}{\lfloor\frac{n}{2}\rfloor !
\lfloor\frac{n}{4}\rfloor ! \lfloor\frac{n}{8}\rfloor !
\cdots}\right) & = & \sum_{i = 1}^{\infty} \left\lfloor
\frac{n}{p^i}\right\rfloor - \sum_{j = 1}^{\infty}\sum_{i =
1}^{\infty}
\left\lfloor\frac{n}{2^j p^i}\right\rfloor  \notag \\
& = & \sum_{i =
1}^{\alpha_p}\left(\left\lfloor\frac{n}{p^i}\right\rfloor -
\sum_{j = 1}^{\alpha_2}\left\lfloor\frac{n}{2^j
p^i}\right\rfloor\right) \label{eq8}
\end{eqnarray}
Next, for all $i \in \{1 , 2 , \dots , \alpha_p\}$, we have:
$$\sum_{j = 1}^{\alpha_2} \left\lfloor\frac{n}{2^j p^i}\right\rfloor ~=~ \sum_{j = 1}^{\alpha_2}
\left\lfloor\frac{\left\lfloor\frac{n}{p^i}\right\rfloor}{2^j}\right\rfloor
~\leq~ \sum_{j = 1}^{\alpha_2}
\frac{\left\lfloor\frac{n}{p^i}\right\rfloor}{2^j} ~<~
\left\lfloor\frac{n}{p^i}\right\rfloor .$$ But since
$(\lfloor\frac{n}{p^i}\rfloor - \sum_{j = 1}^{\alpha_2}
\lfloor\frac{n}{2^j p^i}\rfloor)$ ($i \in \{1 , 2 , \dots ,
\alpha_p\}$) is an integer, it follows that:
$$\left\lfloor\frac{n}{p^i}\right\rfloor - \sum_{j = 1}^{\alpha_2} \left\lfloor\frac{n}{2^j p^i}\right\rfloor ~\geq~ 1
~~~~~~ (\forall i \in \{1 , 2 , \dots , \alpha_p\}) .$$ By
inserting those last inequalities in (\ref{eq8}), we finally
obtain:
$$v_p\left(\frac{n!}{\lfloor\frac{n}{2}\rfloor !
\lfloor\frac{n}{4}\rfloor ! \lfloor\frac{n}{8}\rfloor !
\cdots}\right) \geq \alpha_p ,$$ which confirms (\ref{eq7}) and
completes this proof.\hfill$\blacksquare$

\begin{thm}\label{t3}
For all positive integer $n$, we have:
$$\prod_{r = 1}^{n} h(r) ~\leq~ \c^n ,$$
where $c = 4.01055487\dots$. \\ In addition, the inequality
becomes an equality for $n = 1023 = 2^{10} - 1$.
\end{thm}

\noindent{\bf Proof.} First, we use the relation (\ref{eq6}) to
prove by induction on $n$ that:
\begin{equation}\label{eq9}
\prod_{r = 1}^{n} h(r) ~\leq~ n^{\log_2 n} 4^n
\end{equation}
$\bullet$ For $n = 1$, (\ref{eq9}) is clearly true.\\
$\bullet$ For a given $n \geq 2$, suppose that (\ref{eq9}) is true
for all positive integer $< n$ and let us show that (\ref{eq9}) is
also true for $n$. To do so, we distinguish the two following
cases:\\
{\bf 1\textsuperscript{st} case:} (if $n$ is even, that is $n = 2
m$ for some $m \in \mathbb{N}^*$). \\ In this case, by using
(\ref{eq6}) and the induction hypothesis, we have:
\begin{eqnarray*}
\prod_{r = 1}^{n} h(r) & = & \binom{2 m}{m} \prod_{r = 1}^{m} h(r) \\
& \leq & \binom{2 m}{m} m^{\log_2 m} 4^m \\
& \leq & m^{\log_2 m} 4^{2 m} ~~~~~~~~ \text{(since $\binom{2
m}{m} \leq 4^m$)} \\
& \leq & n^{\log_2 n} 4^n ,
\end{eqnarray*}
as claimed.\\
{\bf 2\textsuperscript{nd} case:} (if $n$ is odd, that is $n = 2 m
+ 1$ for some $m \in \mathbb{N}^*$). \\ By using (\ref{eq6}) and
the induction hypothesis, we have:
\begin{eqnarray*}
\prod_{r = 1}^{n} h(r) & = & (2 m + 1) \binom{2 m}{m} \prod_{r = 1}^{m} h(r) \\
& \leq & (2 m + 1) \binom{2 m}{m} m^{\log_2 m} 4^m \\
& \leq & m^{\log_2 m + 1} 4^{2 m + 1} ~~~~~~~~ \text{(since $2 m +
1 \leq 4 m$ and $\binom{2
m}{m} \leq 4^m$)} \\
& \leq & n^{\log_2 n} 4^n ,
\end{eqnarray*}
as claimed. \\ The inequality (\ref{eq9}) thus holds for all
positive integer $n$. Now, to establish the inequality of the
theorem, we
proceed as follows:\\
--- For $n \leq 70000$, we simply verify the truth of the inequality in
question (by using the Visual Basic language for example).\\
--- For $n > 70000$, it is easy to see that $n^{\log_2 n} \leq
(\c/4)^n$ and by inserting this in (\ref{eq9}), the inequality of
the theorem follows.\\
The proof is complete.\hfill$\blacksquare$~\vspace{2mm}

Now, since any positive integer $n$ satisfies $\prod_{r = 1}^{n}
f(r) \leq \prod_{r = 1}^{n} h(r)$ (according to (\ref{eq5}) and
the fact that $\prod_{r = 1}^{n} h(r)$ is an integer), then we
immediately derive from Theorem \ref{t3} the following:
\begin{coll}\label{coll1}
For all positive integer $n$, we have:
$$\prod_{r = 1}^{n} f(r) ~\leq~ \c^n ,$$
where $\c$ is the constant given in Theorem
\ref{t3}.\hfill$\blacksquare$
\end{coll}

To improve Corollary \ref{coll1}, we propose the following optimal
conjecture which is very probably true but it seems difficult to
prove or disprove it!
\begin{conj}\label{conj1}
For all positive integer $n$, we have:
$$\prod_{r = 1}^{n} f(r) ~<~ 4^n .$$
\end{conj}
Using the Visual Basic language, we have checked the validity of
Conjecture \ref{conj1} up to $n = 100000$. Further, by using
elementary estimations similar to those used in the proof of
Theorem \ref{t3}, we can easily show that:
$$\lim_{n \rightarrow + \infty}\left(\prod_{r = 1}^{n} f(r)\right)^{\!\!1/n} = \lim_{n \rightarrow + \infty}\left(\prod_{r = 1}^{n}
h(r)\right)^{\!\!1/n} = 4 ,$$ which shows in particular that the
upper bound of Conjecture \ref{conj1} is optimal.

Now, by exploiting the properties obtained above for the
arithmetic function $f$, we are going to establish some curious
properties concerning the $2$-adic valuation.

\begin{thm}\label{t4}
For all positive integer $n$ and all odd prime number $p$, we
have:
$$\sum_{r = 1}^{n} v_2(r) v_p(r) ~\leq~ \sum_{r = 1}^{n} v_p(r) - \left\lfloor\frac{\log n}{\log p}\right\rfloor .$$
\end{thm}

\noindent{\bf Proof.} Let $n$ be a positive integer and $p$ be an
odd prime number. Since (according to Theorem \ref{t2}), the
product $\prod_{r = 1}^{n} f(r)$ is a multiple of the positive
integer $\odd(\lcm(1 , 2 , \dots , n))$ whose the $p$-adic
valuation is equal to $\lfloor\frac{\log n}{\log p}\rfloor$, then
we have:
$$v_p\left(\prod_{r = 1}^{n} f(r)\right) = \sum_{r = 1}^{n} v_p\left(f(r)\right) ~\geq~ \left\lfloor\frac{\log n}{\log p}\right\rfloor .$$
But by the definition of $f$, we have for all $r \geq 1$:
$$v_p(f(r)) ~=~ (1 - v_2(r)) v_p(r) .$$
So, it follows that:
$$\sum_{r = 1}^{n} (1 - v_2(r)) v_p(r) \geq \left\lfloor\frac{\log n}{\log p}\right\rfloor ,$$
which gives the inequality of the theorem.\hfill$\blacksquare$

\begin{thm}\label{t5}
Let $n$ be a positive integer and let $a_0 + a_1 2^1 + a_2 2^2 +
\dots + a_s 2^s$ be the representation of $n$ in the binary
system. Then we have:
$$\sum_{r = 1}^{n} \frac{v_2(r) (3 - v_2(r))}{2} ~=~ \sum_{i = 1}^{s} i a_i .$$
In particular, we have for all $m \in \mathbb{N}$:
$$\sum_{r = 1}^{2^m} \frac{v_2(r) (3 - v_2(r))}{2} ~=~ m .$$
\end{thm}

\noindent{\bf Proof.} By taking the 2-adic valuation in the two
hand-sides of the identity (\ref{eq4}) and then using (\ref{eq6}),
we obtain:
$$\sum_{r = 1}^{n} \frac{v_2(r) (3 - v_2(r))}{2} = v_2\left(\prod_{r = 1}^{n} h(r)\right) =
v_2\left(\frac{n!}{\lfloor\frac{n}{2}\rfloor!
\lfloor\frac{n}{4}\rfloor! \lfloor\frac{n}{8}\rfloor!
\cdots}\right) .$$ It follows by using Legendre's formula (see
e.g., \cite{hw}) that:
\begin{eqnarray*}
\sum_{r = 1}^{n} \frac{v_2(r) (3 - v_2(r))}{2} & = & \sum_{i =
1}^{\infty} \left\lfloor\frac{n}{2^i}\right\rfloor - \sum_{j =
1}^{\infty}\sum_{i = 1}^{\infty} \left\lfloor\frac{n}{2^{i +
j}}\right\rfloor \\
& = & \sum_{i = 1}^{\infty} \left\lfloor\frac{n}{2^i}\right\rfloor
- \sum_{u = 2}^{\infty} (u - 1)
\left\lfloor\frac{n}{2^u}\right\rfloor \\
& = & \sum_{i = 1}^{\infty} \left\lfloor\frac{n}{2^i}\right\rfloor
- \sum_{i = 1}^{\infty} i \left\lfloor\frac{n}{2^{i +
1}}\right\rfloor .
\end{eqnarray*}
By adding to the last series the telescopic series $\sum_{i =
1}^{\infty} \left((i - 1) \left\lfloor\frac{n}{2^i}\right\rfloor -
i \left\lfloor\frac{n}{2^{i + 1}}\right\rfloor\right)$ which is
convergent with sum zero, we derive that:
$$\sum_{r = 1}^{n} \frac{v_2(r) (3 - v_2(r))}{2} = \sum_{i = 1}^{\infty} i \left(\left\lfloor\frac{n}{2^i}\right\rfloor
- 2 \left\lfloor\frac{n}{2^{i + 1}}\right\rfloor\right) .$$ But
according to the representation of $n$ in the binary system, we
have:
$$\left\lfloor\frac{n}{2^i}\right\rfloor - 2 \left\lfloor\frac{n}{2^{i + 1}}\right\rfloor
= \begin{cases} a_i & \text{for $i = 1 , 2 , \dots , s$}
\\ 0 & \text{for $i > s$}\end{cases} .$$
Hence
$$\sum_{r = 1}^{n} \frac{v_2(r) (3 - v_2(r))}{2} = \sum_{i = 1}^{s} i a_i ,$$
as required.\\
The second part of the theorem is nothing else an immediate
application of its first part with $n = 2^m$. The proof is
finished.\hfill$\blacksquare$

\end{document}